\documentstyle[12pt,epsf]{article}

\title{The index of $grad \, f(x,y)$
\\(Revised version
\footnote{A thoroughly revised and hopefully more readable version of Duke eprint alg-geom
9506002; the main results are the same.} )}
\author{Alan H. Durfee}

\newcommand{\postscript}[2]
{\setlength{\epsfxsize}{#2\hsize}
\centerline{\epsfbox{#1}}}

\def\endofproof{$\Box$}

\def\real{{\bf R}}
\def\complex{{\bf C}}
\def\projective{{\bf P}}

\def\complexlineinfinity{{\bf L}_\complex}
\def\clineinf{\complexlineinfinity}

\def\reallineinfinity{{\bf L}}
\def\rlinf{\reallineinfinity}
\def\lineinf{\reallineinfinity}

\def\circleinf{{\bf S}}

\def\infinity{\infty}

\def\euler{\chi}

\def\intersect{\cdot}

\def\tf{\tilde{f}}
\def\realinf{\real \cup \{\infinity \} }
\def\dr{d_{{\bf R}}}
\def\flevels{f \cap \reallineinfinity}

\def\drt{\tilde{d}_{{\bf R}}}

\newcounter{mycounter}[section]
\renewcommand{\themycounter}{\arabic{section}.\arabic{mycounter}}

\newenvironment{theorem}%
{\medskip 
    \refstepcounter{mycounter}
    {\bf \noindent Theorem \themycounter.}\begin{em}}%
{\end{em} \medskip }

\newenvironment{proposition}%
{\medskip 
    \refstepcounter{mycounter}
    {\bf \noindent Proposition \themycounter. \ } \begin{em} }%
{\end{em} \medskip }

\newenvironment{lemma}%
{\medskip 
    \refstepcounter{mycounter}
    {\bf \noindent Lemma \themycounter. \ } \begin{em} }%
{\end{em} \medskip }

\newenvironment{corollary}%
{\medskip 
    \refstepcounter{mycounter}
    {\bf \noindent Corollary \themycounter. \ } \begin{em} }%
{\end{em} \medskip }

{\medskip 
    \refstepcounter{mycounter}
    {\bf \noindent Formula \themycounter. \ } \begin{em} }%
{\end{em} \medskip }

{\medskip 
    \refstepcounter{mycounter}
    {\bf \noindent Remark \themycounter. \ }}%
{\medskip }

\newenvironment{conjecture}%
{\medskip 
    \refstepcounter{mycounter}
    {\bf \noindent Conjecture \themycounter. \ }}%
{\medskip }

\newenvironment{definition}%
{\medskip 
    \refstepcounter{mycounter}
    {\bf \noindent Definition \themycounter. \ }}%
{\medskip }

\newenvironment{example}%
{\medskip 
    \refstepcounter{mycounter}
    {\bf \noindent Example \themycounter. \ }}%
{\medskip }

{\medskip 
    \refstepcounter{mycounter}
    {\bf \noindent Problem \themycounter. \ }}%
{\medskip }

\newenvironment{xproof}%
{\medskip 
  \noindent
    {\bf Proof. \ }}%
{\endofproof \medskip }

\begin{document}
\maketitle

\begin{abstract}
Let $f(x,y)$ be a real polynomial of degree $d$ with isolated critical
points, 
and let $i$ be the index of $grad \, f$ around a large circle
containing the critical points.
An elementary argument shows that $|i| \leq d-1$.
In this paper we show that $ i \leq max \{ 1, d-3 \}$.
We also show that
if all the level sets of $f$
are compact, then $i = 1$, and otherwise
$ |i| \leq \dr -1 $
where  $\dr$ is the sum of the multiplicities of the real linear
factors in the homogeneous term of highest degree in $f$.  
The technique of proof involves computing $i$ from information at
infinity.
The index $i$ is broken up into a sum of components
$i_{p,c}$ corresponding to points $p$ in the real 
line at infinity and
limiting values $c \in \realinf$ of
the polynomial. 
The numbers $i_{p,c}$ are computed in three ways:  geometrically, from a
resolution of $f(x,y)$, and from a Morsification of $f(x,y)$.
The $i_{p,c}$ also provide
a lower bound for the number of
vanishing cycles of $f(x,y)$ at the point $p$ and value $c$.
\end{abstract}



\section{Introduction}

Let $f(x,y)$ be a real polynomial with isolated critical points.
Let $i$ be the index of the gradient vector field of $f(x,y)$ around a
large circle $C$ centered at the origin and containing the critical
points, oriented in the counterclockwise direction.  If the critical
points of $f$ are nondegenerate, then the index $i$ is the number of
local extrema minus the number of saddles.

What bounds can be placed on the index $i$ in terms of the degree $d$
of the polynomial?
It follows easily from Bezout's theorem that \cite[Proposition 2.5]{REU}
$$\vert i \vert \leq d-1$$
 
It is easy to find polynomials satisfying the lower bound of this
inequality; for example if $f = l_1 \dots l_d$ where the $l_i$ are
equations of lines in general position, then $i = 1-d$, as can be seen
by looking at how the gradient vector field turns on the circle $C$,
or by counting critical points \cite[Section 4]{REU}.

The upper bound is more mysterious.
In the first place, polynomials with $i >1$ are hard to find. 
(The dubious reader should try to do so!)
A simple example with two local extrema and no other critical points 
($i=2$) is
$f(x,y) = y^5+x^2y^3-y$.
A polynomial of degree five can have as many as sixteen critical
points in the complex plane; a generic polynomial of degree five will
have exactly this number.
The above polynomial, however, has only
four critical points in the plane (two real and two complex), so it is
not generic.
In fact this behavior is typical for polynomials with $i>1$
\cite[Theorem 6.2]{REU}.

There are polynomials of degree $d$ with $i$ arbitrarily large
(see Example \ref{many-max-min-ex}), but they have $i \approx (1/3)d$.
So evidently there is a large gap between the theoretical upper bound
and examples.
One of the goals of this paper is to give a modest improvement of this
upper bound.  We will show

\medskip

\noindent
{\bf Theorem \ref{max-theorem}.}  If $f(x,y)$ is a real polynomial of degree $d$ with isolated
critical points, and $i$ is the index of $grad \, f$ around a large
circle containing the critical points, then
$$ i \leq max \{ 1, d-3 \}$$

\medskip

In particular this result implies that the minimum degree for a polynomial
with $i > 1$ is five, as in the example above.
In fact, the bound is often better. (See, for example, Proposition \ref{i-upper-bound}.)

Let the {\em real degree} $\dr$ of $f$ be the sum of the multiplicities of the real linear
factors in the homogeneous term of highest degree in $f$.
(Thus $\dr \leq d$.)  
We will also show 

\medskip

\noindent
{\bf Theorem \ref{dr-theorem}.}  If all the level sets of the polynomial $f(x,y)$
are compact, then $i = 1$.  Otherwise
$$|i| \leq \dr -1 $$

\medskip

It is easy to find polynomials realizing the lower bound (Corollary
\ref{corollary-dr}), but the upper bound still appears high.

The basic idea of the proofs is to compute the index $i$ from ``information at
infinity''.  We write $i$ as 
$$i = 1 + \sum_{\scriptstyle p \in \rlinf \atop \scriptstyle c \in \realinf} i_{p,c}$$
The terms $i_{p,c}$ are defined as follows:
The number $\pm 1/2$ is assigned to a point $q$ where the circle $C$ is tangent
to a level set of the polynomial according as whether the level set is
locally inside or outside $C$ at $q$.  The circle is then made larger and
larger.  The point $q$ where the level set is tangent to the circle approaches a limiting point $p$ on the
line at infinity in real projective space, and the value of the polynomial $f(q)$ approaches a
limiting value $c$.  
The term $i_{p,c}$ is the sum of all the numbers $\pm 1/2$ associated
to $p$ and $c$ in this manner.
This material is in Section 3.
We also show (Proposition \ref{generic-h}) that the family of circles can be
replaced by the level sets of any reasonable function, and the
$i_{p,c}$ will remain the same.

The polynomial $f$ extends to a function on projective space which is
not well-defined at certain points on the line at infinity.
Blowing up these points gives a well-defined function $\tf$.
We use this technique in Section 4 to derive some simple properties of the level curves of
$f$.  

In Section 5, we use Morse theory to show that the $i_{p,c}$ can be computed from
the critical points of $\tf$ and information about the exceptional
sets.
The process of blowing up and computing the index is easy to carry out
in specific examples.

The polynomial $f$ can also be deformed into what we call a
``Morsification'', a polynomial whose real critical points are
nondegenerate and whose homogeneous term of highest degree has no
repeated real linear factors (Section 6).  There is a simple formula
relating the index of the original polynomial, the index of the new
polynomial, and the index of the newly created critical points.  The
deformation process is not too well understood, and this section
contains some examples and a conjecture.

The computations of these sections are used in Section 7 to establish bounds on
the $i_{p,c}$.
These local bounds are sharp.
The global bounds on $i$ follow from the local bounds and some 
delicate arguments.
However, the global bounds are not sharp and
there still is a big gap between the global bounds and the
examples.

In Section 8 we relate $i_{p,c}$ to the ``jump'' at $c$ in the Milnor
number of the family $f(x,y) = t$  at the point $p$ on the
line at infinity.

Throughout this paper the techniques are those of basic topology
(Morse Theory) and basic algebraic geometry (Bezout's theorem,
explicit computation of intersection multiplicities, etc.)  Computer
algebra programs were used to find critical points, countour plots and
the $i_{p,c}$.  Although many of the results and techniques are valid
in higher dimensions, the exposition is in dimension two for reasons
of clarity.

The author's interest in these questions started in 1989 when he worked with a
group of undergraduates in the Mount Holyoke Summer REU \cite{REU}.
Another group of students continued this work in 1992; one of their
results was the construction of polynomials with an arbitrarily large
number of local maxima and no other critical points
\cite{Robertson-P}.  (These polynomials have $ i \approx d/4$.)

Shustin \cite{Shustin-96} has studied polynomials all of whose critical
points lie in the complex plane.   
He finds polynomials of this type with almost all arbitrarily prescribed numbers of
local maxima, minima and saddles.
These polynomials have $i = 1-\dr$ and, in particular, $i \leq 1$.
They are stable in the sense that nearby polynomials have the same
number and type of critical points.  The primary focus of this paper
is polynomials $f$ with $i > 1$; these polynomials
are not stable.
In fact, \cite[Theorem 6.2]{REU} says that
$$i \leq \frac{1}{2} m + 1$$ 
where $m$ is the sum over $p$ in the line at infinity in real
projective space of the intersection multiplicities at $p$ of the
completions of $f_x = 0$ and $f_y = 0$.

This paper is real counterpart of the study by many people of
``critical points at infinity'' for complex polynomials; see
\cite{Durfee-P96} for further references.

Research on this paper was partially supported by NSF grant
DMS-8901903, and a grant from the
International Research and Exchanges Board (IREX), with funds provided
by the National Endowment for the Humanities and the United States
Information Agency.  
The research was carried out over the past eight years 
at Martin-Luther University, Halle,
the University of Nijmegen,
Warwick University and the Massachusetts Institute of Technology; the
author would like to thank them for their hospitality.

Some notation which will be used throughout the paper:
We let 
$$\reallineinfinity = \{ [x,y,z] \in \projective^2: z=0 \}$$ 
be the line at infinity in real projective space $\projective^2$,
and $\clineinf$ be the line at infinity in complex projective space
$\complex \projective ^2$.
We use $d$ for the degree of the polynomial $f(x,y)$, and
$f_d$ for the homogeneous term of degree $d$ in $f$.



\section{A polynomial zoo}

A number of polynomials with strange properties are used as examples throughout
this paper.  These are described in this section.

\begin{example}
\label{std-crpt-ex}
The polynomial  $y(xy-1)$, which has no critical points in the plane, 
is the standard example of a polynomial with a
``critical point at infinity'' (at $[1,0,0]$).
The ``critical value'' (jump in the Milnor number) is at 0.
This polynomial perhaps first appeared in \cite{Broughton-83}.
\end{example}

\begin{example}
\label{std-nocrpt-ex}
The polynomial $x(y^2-1)$ has saddles at $(0,1)$ and $(0,-1)$.
The family of level curves at $[1,0,0]$ is equisingular;
there is no ``critical point'' at $[1,0,0]$.
\end{example}

\begin{example}
\label{parabola-ex}
The parabola $y^2-x$ is the simplest example of a polynomial with a
``critical point'' at $[1,0,0]$ with ``critical value'' $\infinity$
\cite{Durfee-P96}.
\end{example}

\begin{example}
\label{max-min-ex}
The polynomial $y^5 + x^2y^3 - y$ has a local minimum at
$(0,-1/\sqrt[4]{5})$, a local maximum at  $(0,1/\sqrt[4]{5})$
and no other critical points.  
This polynomial was found by an REU group of undergraduates at Mount Holyoke
College in the summer of 1989 \cite{REU}.
\end{example}

\begin{example}
\label{many-max-min-ex}
The polynomial $(y(x^2+1)-1)(y(x^2+2)-1) \dots (y(x^2+k)-1)$ 
has $k-1$ local extrema
and no other real critical points; for $k = 2$ there
is a local minimum, for $k=3$ there is a local minimum and a local
maximum, for $k=4$ there are two local minima and a local maximum, and
so forth. This polynomial was also found by the REU group \cite{REU}.
\end{example}

\begin{example}
\label{two-min-ex}
The polynomial $(xy^2 - y -1)^2 + (y^2-1)^2$ from \cite{Mathmag-85} has
local minima at $(2,1)$ and $(0, -1)$, and no other critical
points.  Note the asymmetry of this polynomial compared with the
previous ones.  
\end{example}

\begin{example}
\label{only-crpt-ex}
The polynomial $x^2(1+y)^3+y^2$ has its sole critical point at the
origin.  
This critical point is a local maximum, but not an absolute
maximum \cite{Calvert-Vam}.
\end{example}

\begin{example}
\label{kras-ex}
The polynomial $y -(xy-1)^2$ has a saddle at
$(-1/2,0)$ and no other critical points.  At $[1,0,0]$ the level set
$f = 0$ has one branch, but the general level set has two branches \cite{Krasinski-91}.
\end{example}

\begin{example}
\label{two-parabola-ex}
The polynomial $f(x,y) = (x-y^2)((x-y^2)(y^2+1) - 1)$ has its zero
locus along the parabola $x = y^2$ and the curve $x = y^2 + 1/(y^2+1)$
which is asymptotic to this parabola.  Its only critical point is a
minimum at $(1/2,0)$.  The level curves intersect $\rlinf$
only at $[1,0,0]$, and they are tangent to $\rlinf$ at this point.
(The ``curve of tangencies'' (see the next section) is also tangent to
$\rlinf$ at $[1,0,0]$.) 
\end{example}


\section{A formula for $i$ from the geometry of $grad \, f$}

Let $f(x,y)$ be a real polynomial with isolated critical points.
(Note that $f$ is thus not constant.)
Let $i$ be the index of the gradient vector field of $f(x,y)$ around a
large circle $C$ centered at the origin and containing the critical points, oriented in the
counterclockwise direction.  (Recall that the index is the topological
degree of the map $C \to S^1$ defined by $ t \mapsto grad \,
f(\alpha(t)) / |grad \, f(\alpha(t))|$, where $t \mapsto \alpha(t) $ is
a parameterization of $C$.)  
This section contains the fundamental geometric decomposition
of the index $i$ (Proposition \ref{geometric-index-formula}).

\begin{figure}
\postscript{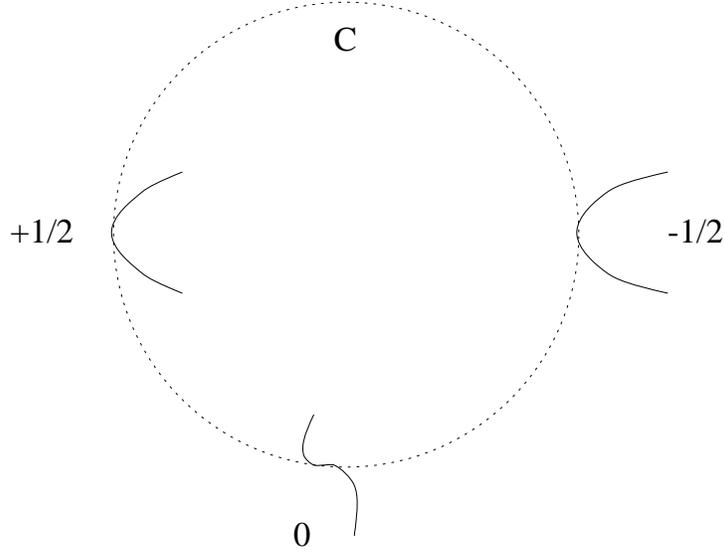}{0.7}
\caption{Assigning $\pm 1/2$ or 0 to a point of tangency}
\label{u-sign}
\end{figure}

To each point $q \in C$ where a smooth level curve of $f$ is
tangent to $C$ at $q$ 
we assign the number $\pm 1/2$ or 0 as follows:
If the level curve of $f$ is
outside the circle $C$ near $q$, this number is $-1/2$.
If it is inside $C$ near $q$, the number is $+1/2$.
(These conditions are topological; the tangency may be algebraically degenerate.)
If one side is outside and the other inside, or if the level set is
contained in $C$ near $q$ (in which case $C$ is a
connected component of the level set), the number is $0$.
(See Figure \ref{u-sign}; the circle $C$ is dotted, and the level
curves of $f$ are solid lines.)

The points $q$ where the level sets of $f$ are tangent to $C$ are the
zeros of the (real) {\em curve of tangencies} 
$$f_x y - f_y x = 0$$
This curve may have reducible components.

Choose the circle $C$ large enough so that it 
contains the compact components and the isolated singular points of the curve
of tangencies and their points of common tangency.
In the exterior of $C$ the curve of tangencies is a union of
connected components.
Each component $\gamma$ is a smooth arc which goes to infinity; we
call this an {\em end} of the curve of tangencies.
Choose the circle $C$ large enough so that the numbers $\pm 1/2, 0$
assigned above are constant along each end $\gamma$.  (This is possible
since the intersection multiplicity of $C$ and the level sets of $f$
is constant along each end $\gamma$ for $C$ large.)
Let $i(\gamma)$ be the number $\pm 1/2$ or $0$ assigned to $\gamma$ in this fashion.

Let $p(\gamma) \in \lineinf$ be the endpoint of the closure of
$\gamma$, and 
let $c(\gamma) \in \realinf$ be the limiting value of $f(q)$ as $q$
goes to infinity along $\gamma$.  

\begin{lemma}
\label{lemma-cgamma}
For each end $\gamma$ of the curve of tangencies, the number
$c(\gamma)$ exists.  In fact, the function $f$ restricted to $\gamma$
is strictly increasing or decreasing.
\end{lemma}

\begin{xproof}
Let
$\Gamma(f) \subset \projective^2 \times \projective$ be the
closure of the graph of $f$.
The end $\gamma$ lifts uniquely to $\Gamma(f)$,
intersecting the fiber over $p$ at a point $(p,c)$.
The number $c$ is $c(\gamma)$.
The function $f$ is strictly increasing or decreasing since $\gamma$ is
perpendicular to the level sets of $f$.
\end{xproof}

For $p \in \reallineinfinity$ and $c \in \realinf$, we let
$$i_{p,c} = \sum i(\gamma)$$
where the sum is over all ends $\gamma$ with $p(\gamma) = p$ and $c(\gamma)
= c$.
We also let
$$i_p = \sum_{c \in \realinf} i_{p,c}$$
and
$$i_{\reallineinfinity, \infinity} = \sum_{p \in \rlinf} i_{p,\infinity}$$

\begin{lemma}
The numbers $i_{p,c}$ are integers (not just half-integers).
\end{lemma}

\begin{xproof}
The curve of tangencies can be lifted to $\Gamma(f)$,
the graph of $f$.
A real branch of this curve at $(p,c) \in \Gamma(f)$ is a pair of ends $\gamma \neq \gamma'$.
If $i(\gamma) = 0$, then the intersection of the level sets of $f$
with the family of circles is degenerate along $\gamma$, and hence
$i(\gamma') = 0$ as well.
\end{xproof}

\begin{proposition}
\label{geometric-index-formula}
If $f(x,y)$ is a real polynomial with isolated critical points, then
$$i = 1 + \sum_{\scriptstyle p \in \rlinf \atop \scriptstyle 
c \in \real} i_{p,c} + i_{\reallineinfinity, \infinity}$$
\end{proposition}

\begin{xproof}
If no connected component of a level set of $f$ is contained in the
large circle $C$, then we have that 
\begin{equation}
\label{i-eqn}
i = 1 + \sum i(\gamma) 
\end{equation}
where the sum is over all ends $\gamma$ of the curve of tangencies:
This is clearly true if all the points of tangency are regular values
for the map $C \to S^1$ defined by $ t \mapsto grad \,
f(\alpha(t)) / |grad \, f(\alpha(t))|$, where $t \mapsto \alpha(t) $ is
a parameterization of $C$.
If a point of tangency is not a regular value for this map (eg for
$f(x,y) = y^3 + x$, or $y^4 + x$), then a small (topological)
deformation shows that it still holds.
The expression of the proposition is just a decomposition of (\ref{i-eqn}). 

Now suppose that a connected component 
of a level set $f=c$ is contained in $C$.
We may assume without loss of generality that
$c \gg 0$.
Since one component of $f=c$ is compact, 
all components are compact by Proposition \ref{geom-prop-noncpt}.
Thus by Proposition \ref{index-prop}, $i = 1$ (which is
obvious here),
$i_{\reallineinfinity, \infinity} = 0$ and $i_{p,c} = 0$ for all $p
\in \reallineinfinity$ and $c \in \real$. 
\end{xproof}

A corollary of the Proposition is that
\begin{equation}
\label{i-ip-formula}
i = 1 + \sum_{p \in \rlinf} i_p
\end{equation}

\begin{figure}
\postscript{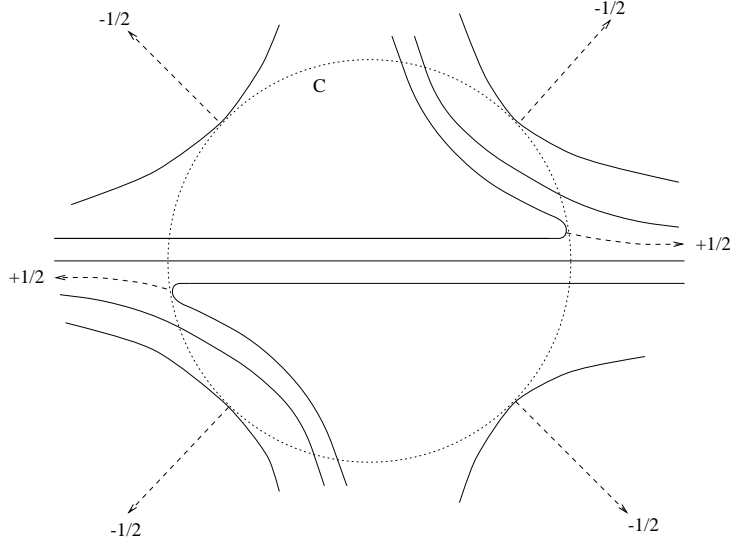}{0.7}
\caption{The index computation for the polynomial $y(xy-1)$}
\label{std-crpt-index}
\end{figure}

The process of decomposing the index for the polynomial $f(x,y) =
y(xy-1)$ of Example \ref{std-crpt-ex} is
pictured in Figure \ref{std-crpt-index}.  (The circle is dotted, the solid
lines are the level sets of $f$, the dashed lines are the
ends $\gamma$ of the curve of tangencies, and the numbers are
$i(\gamma)$.)

A geometrically obvious example of the decomposition of the
index is for a polynomial $f$ for which the real linear factors of 
$f_d$ are irreducible.  In this
case $i = 1 - \dr$, $i_{p,c}= 0$ for $p \in \rlinf$ and $c \in \real$,
and $i_{\rlinf, \infinity} = -\dr$.  (This will be proved formally in
Corollary \ref{corollary-dr}.)

The invariants of Proposition \ref{geometric-index-formula} for selected
polynomials 
are given in Table 1; all the nonzero $i_{p,c}$ for $c \in \real$ are listed.

\begin{table}
\label{index-invariants}
\begin{tabular}{|l|r|l|l|l|r|}  \hline
$f(x,y)$ & $i_{\reallineinfinity, \infinity}$ & $p \in \flevels$ & $c$ & $i_{p,c}$  & $i$  \\ \hline \hline 
Example \ref{std-crpt-ex}: $y(xy-1)$ &$-2$ & $[1,0,0]$ & $0$ & 1 & 0  \\ \hline 
Example \ref{std-nocrpt-ex}: $x(y^2-1)$ &$ -3$ & & & &$-2$ \\ \hline
Example \ref{parabola-ex}:  $y^2-x$ &$-1$ & & & & 0 \\ \hline
Example \ref{max-min-ex}:  $y^5 + x^2y^3-y$ & $-1$ & $[1,0,0]$ & $0$ & 2 &2 \\ \hline 
Example \ref{two-parabola-ex} & $-1$ & $[1,0,0]$ & $0$ & 1 & 1 \\ \hline
Example \ref{kras-ex}:  $y - (xy-1)^2$ & $-2$ & $[1,0,0]$ & 0 & 0 &  $-1$ \\  \hline 
Example \ref{two-min-ex} & $-1$ & $[1,0,0]$ & 1 & 1 & 2 \\ \cline{4-5}
                 &   &   & 2 & 1 &   \\ \hline 
$y(x^2y-1)$ & $-2$ & $[1,0,0]$ & 0 & 1 & 0 \\ \hline
\end{tabular}
\caption{Index invariants of selected polynomials}
\end{table}

Note that the sum of the $i(\gamma)$'s making up $i_{p,c}$ is over
ends $\gamma$ 
where $grad \, f$ points both out of and into the circle $C$;
the process of decomposing the index described below does not work if
the sum is just over those points where the gradient points out, as
can be seen in the example $f(x,y) = y(x^2y-1)$.

It is useful to have both the expression of Proposition
\ref{geometric-index-formula} where the limiting value $c=\infinity$
is separated out and put into
$i_{\lineinf, \infinity}$ (see, for example, Proposition
\ref{res-index-formula}), as well as the expression of Equation
(\ref{i-ip-formula}), where these values are grouped by $p$ into $i_p$
(see, for example, Lemma \ref{i-d-local-estimate}).

The decomposition of $i$ into the $i_{p,c}$ reflects the geometry of
$f$ near infinity.  It is apparently not related to the finite
critical points of the polynomial and their critical values.



\section{Resolutions and the geometry of level sets}

In this section we describe the resolution of the points of
indeterminacy of a polynomial on the line at infinity, and use this
concept to establish some simple properties of its affine level curves. 

A polynomial
$$f: \real^2 \to \real $$ extends to a map of real projective spaces 
$$\hat{f} : \projective^2 \to \projective $$ which is undefined at a
finite number of points on the line at infinity $\reallineinfinity$.
By blowing up these points one gets a manifold $M$ and a map $$ \pi : M \to \projective^2 $$ such
that the map $$ \tilde{f} : M \to \projective $$ lifting
$\hat{f}$ is everywhere defined.  We call the map $ \tilde{f}$ a
{\em resolution of $f$}.  
(We avoid the use of minimal resolutions, though this concept could be
used to
provide alternate proofs of some of the results below.)



Any resolution $ \tilde{f} : M \to \projective $ factors through 
$$\bar{f} : \Gamma(f) \to \projective $$
where $\Gamma(f)$ is the graph of $f$ as defined above; 
note that the function $\bar{f}$ is everywhere defined.
 
\begin{figure}
\postscript{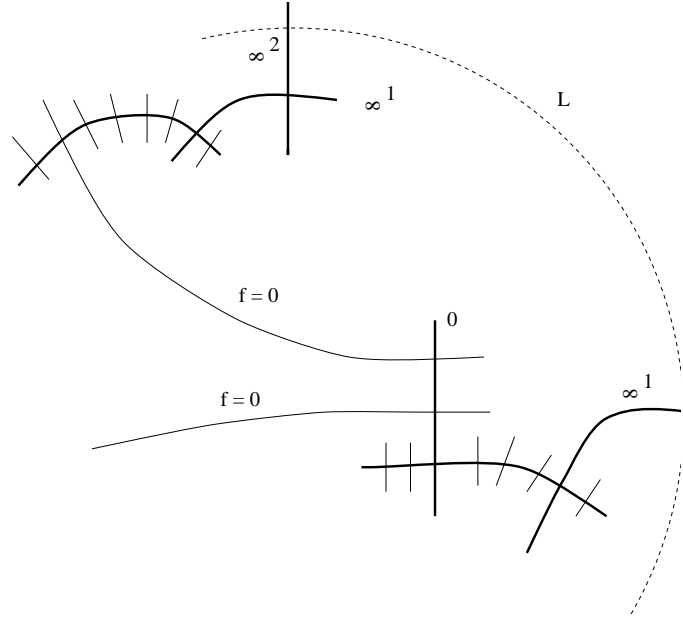}{0.7}
\caption{A resolution of $y(xy-1)$}
\label{std-crpt-res-g}
\end{figure}

For example, a resolution of 
$y(xy-1)$ is given in Figure \ref{std-crpt-res-g}.
The proper transforms of level curves of $f$ are ordinary lines;
the exceptional sets are thick lines.
The number $c^m$ next to a divisor means that at each smooth
point of the divisor there are local coordinates $(u,v)$ in a
neighborhood of the point such that the divisor is $u=0$ and
$\tilde{f}(u,v) =(u-c)^m$. 

Let $f$ be a polynomial, and let $\tf$ be a resolution of $f$.
By $A \gg 0$, we mean as usual that $A$ is large, but more precisely
in this context we mean that $A$ is greater than the absolute value
of all the critical values of $f$, and that if 
$|t| \geq A$, then the
level sets $\tf = t$ are smooth and transversally intersect the
exceptional sets of $\tf$.
In particular, this means that the topological
type of the level sets $f(x,y) = A$ and $f(x,y) = -A$ are independent of $A$.

Fix a polynomial $f$.  Let
$$\lineinf' = \{ [a,b,0] \in \lineinf : f_d(a,b) = 0 \} $$
and let
$$\lineinf'' = \{ p \in \lineinf : \mbox{ There is a $t \in \real$ with
$p$ in the closure of $f(x,y) = t$} \} $$
Both $\lineinf'$ and $\lineinf''$ are finite sets of points.

\begin{lemma}
\label{geom-lemma-1}
If $p   \in \lineinf''$, then 
\begin{enumerate}
\item The point $p$ is in the closure of
level sets $|f| = A$ for $A \gg 0$.
\item If $\tf$ is a resolution of $f$, there is an exceptional set
over $p$ on which $\tf$ is not constant.
\item $p \in \lineinf'$.
\end{enumerate} 
\end{lemma}

\begin{xproof}
Choose a resolution $\tilde{f}$ of $f$.  
There is an exceptional set $E$ with $\pi(E) = p$ such that $\tilde{f}
= t$ intersects $E$.
The function $\tilde{f}$ restricted to $E$ is a real rational
function.

If $\tilde{f}$ restricted to $E$ is not constant, then there is a $q \in E$
such that $\tilde{f}$ restricted to $E - \{q\}$ is a polynomial.
This implies (1) and (2) in this case.

Now suppose that $\tilde{f}$ restricted to $E$ is constant.
The (real) exceptional sets over $p$ form a connected tree.
Since the value of $\tilde{f}$ is $t$ at one point on the tree, and is
infinity at the points where the tree intersects the proper
transform of $\reallineinfinity$, there is a component in the tree
where $\tilde{f}$ takes a continum of large values. This implies (1)
and (2)
in this case.

Part (3) follows since the level curves in $\complex^2$ of the
complexified polynomial $f$ intersect $\clineinf$ at exactly the zeros of
the complexified $f_d$.
\end{xproof}
 
Note that the converse of (3) above is not true:  For the polynomial
$f(x,y) = y^4 + x^2$, for example, 
$f_4(1,0) = 0$, but $[1,0,0]$ is not contained in the
closure of any real level curve.


Given a real polynomial $f(x,y)$, we let $l''$ be the number of points in $\lineinf''$.
If $\tf$ is a resolution of $f$, we let $\xi_{\reallineinfinity,nc}(\tilde{f})$ be the number of (real)
exceptional sets $E$ of $\tf$ such that $\tilde{f} | E$ is nonconstant.
(In the pictures, these exceptional sets are cross hatched by level
curves of the polynomial.)


\begin{corollary}
\label{l-xi}
If $\tf$ is a resolution of $f$, then $l'' \leq \xi_{\rlinf, nc}(\tf)$.
\end{corollary}


The inequality may be strict:
The polynomial 
$x(y+1)(y+2) \dots (y+k)$ 
has $l''=2$ and $\xi_{\rlinf, nc}(\tf) = k+1$.

\begin{proposition}
\label{geom-prop-cpt}
Let $f(x,y)$ be a real polynomial, let $\tf$ be a resolution of $f$,
and let $A \gg 0$.
The following are equivalent:
\begin{enumerate}
\item The set $|f| = A$ is compact.
\item The set $f = t$ is compact, for all $t$.
\item $\xi_{\rlinf, nc}(\tf) = 0$.
\end{enumerate}
Furthermore, if any of the above is true, then
\begin{enumerate}
\item One of $f = A$ and $f = -A$ is homeomorphic to a
circle, and the other is empty.
\item $f_d$ has no real linear factors ($\dr = 0$).
\end{enumerate}
\end{proposition}

\begin{xproof}
(1) implies (2) by Lemma \ref{geom-lemma-1}.
(2) is equivalent to (3) since $f=t$ is compact for all $t$ if
and only if $\tilde{f} | E =
\infinity$ for all exceptional sets $E$.
The additional conclusions are obvious.
\end{xproof}

\begin{proposition}
\label{geom-prop-noncpt}
Let $f(x,y)$ be a real polynomial, let $\tf$ be a resolution of $f$,
and let $A \gg 0$.
If $|f| = A$ is not compact, then
\begin{enumerate}
\item All the connected components of $|f| = A$ are noncompact.
\item The number of connected components of $|f| = A$ is $2\xi_{\rlinf, nc}(\tf)$.
\end{enumerate}
\end{proposition}

\begin{xproof}
The geometry of the resolution
implies (1).
If $E$ is an exceptional set on which $\tilde{f}$ is not constant,
then $\tilde{f}$ restricted to $E$ either takes the value $+A$ exactly
twice, the value $-A$ exactly twice or the values $+A$ and $-A$ each
once.  
Since the level sets $f = \pm A$ are transverse to $E$, this proves
(2).
\end{xproof}


\section{A formula for $i$ in terms of a resolution}

This section gives a formula (Proposition \ref{res-index-formula}) for
computing the index $i$ and the terms $i_{p,c}$ in the decomposition
of $i$ in terms of a
resolution of the polynomial.  In Section 7 this proposition
will
play a  role in finding bounds on $i$.

\begin{lemma}
\label{lemma-i-p-c}
Let $c \in \real$.
If $p \in \reallineinfinity$ is not in the closure of the level set $f
= c$, 
then $i_{p,c} = 0$.
(Furthermore,
$i^{abs}_{p,c} = 0$, in the notation of Section 7.)
\end{lemma}

\begin{xproof}
Suppose $i_{p,c} \neq 0$.
There is an end $\gamma$ of the curve of tangencies passing through $p$, and $f$ has
limiting value $c$ along $\gamma$.
Let $\tilde{f}: M \to \projective^2$ be a resolution of $f$.
The curve $\gamma$ lifts to $M$ and passes through some $q \in
M$ with $\pi(q) = p$.
Also $\tilde{f}(q) = c$.
Hence the closure of $f=c$ intersects some exceptional set over $p$,
so $p$ is in the closure of $f=c$.
This is a contradiction.
\end{xproof}

The precise meaning of the notation $A \gg 0$ can be found in the
previous section.

\begin{proposition}
\label{index-prop}
If $f(x,y)$ is a real polynomial with isolated critical points, and if 
$|f| = A$ is compact for $A \gg 0$, then
\begin{enumerate}
\item $i = 1$.
\item $i_{\reallineinfinity, \infinity} = 0$.
\item $i_{p,c} = 0$ for all $p \in \reallineinfinity$ and $c \in \real$.
\end{enumerate}
\end{proposition}

\begin{xproof}
By Proposition \ref{geom-prop-cpt}, either $f(x,y) = A$ or $f(x,y) = -A$ is 
homeomorphic to a circle; let us assume the former.
Clearly $i = 1$ and $i_{\reallineinfinity, \infinity} = 0$.
Also $f(x,y) = c$ is compact for all $c$, so \ref{lemma-i-p-c} implies
that $i_{p,c} = 0$.
\end{xproof}

Let $\tf$ be a resolution of $f$.
For $p \in \reallineinfinity$ and $c \in \real \cup \{\infinity\}$, let
\begin{itemize}
\item $i_{p,c}(\tilde{f})$ be the sum of the indices of $\tilde{f}$ at
critical points $q \in M$ of $\tf$ such that $\tilde{f}(q) = c$ and $\pi(q) = p$.
\item $\xi_{p,c}(\tilde{f})$ be the number of (real) exceptional sets $E$ of
$\tf$ with $\pi(E) =p$ and $\tilde{f}|E = c$. 
\end{itemize}

Recall that $\xi_{\reallineinfinity,nc}(\tilde{f})$ is the number of (real)
exceptional sets on which $\tilde{f}$ is nonconstant.

\begin{proposition}
\label{res-index-formula}
If $f(x,y)$ is a real polynomial with isolated critical points and
$\tf$ is a resolution of $f$, then
\begin{enumerate}
\item $i = 1 - \sum_{\scriptstyle p \in \rlinf \atop \scriptstyle c
\in \real} \bigl( i_{p,c}(\tf) + \xi_{p,c}(\tf) \bigr) - \xi_{\reallineinfinity,nc}(\tilde{f})$
\item $ i_{p,c} = - i_{p,c}(\tilde{f}) - \xi_{p,c}(\tilde{f})$, for
$p \in \reallineinfinity$ and $c \in \real$.
\item $ i_{\rlinf,\infinity}  = - \xi_{\rlinf,nc}(\tilde{f})$
\end{enumerate}
\end{proposition}

By Proposition \ref{geometric-index-formula}, any two parts of this proposition imply the third, but proving each
part separately is more instructive.
In fact, Part (1) follows from a straight-forward application of Morse theory, Part (2)
follows from Morse theory on a manifold with boundary,
and Part (3) follows from the geometry of the large level curves.

\begin{figure}
\postscript{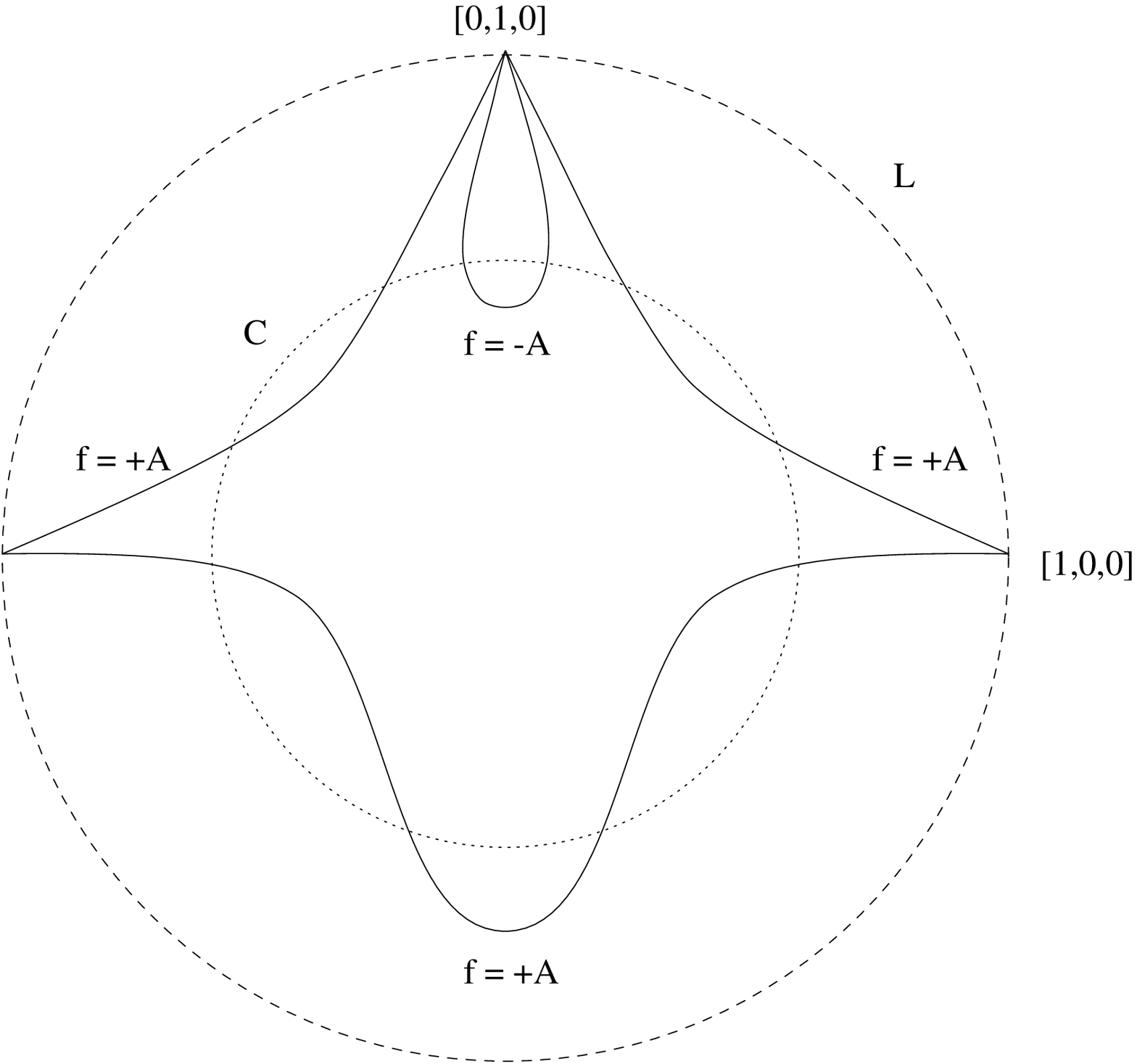}{0.7}
\caption{The level curves of $f(x,y) = y(x^2y-1)$ in $\projective^2$}
\label{large-level-curves}
\end{figure}

\begin{xproof}
Proof of (1):
Let
\begin{itemize} 
\item $i_{\reallineinfinity}(\tilde{f}) = 
\sum_{\scriptstyle p \in \rlinf \atop \scriptstyle c \in \real} i_{p,c}(\tf)$
\item $\xi_{\reallineinfinity,\infinity}(\tilde{f}) = 
\sum_{p \in \rlinf} \xi_{p,\infinity}(\tf)$
\item $\xi_{\rlinf}(\tilde{f})= \sum_{\scriptstyle p \in \rlinf \atop 
\scriptstyle c \in \real} \xi_{p,c}(\tf)$
\end{itemize}
We will do Morse theory on the function $\tf: M \to \real \cup \{
\infinity \} $.  Suppose $ c_1 < c_2 \dots < c_r $ in $\real$ are the
critical values of $\tf$ restricted to the inverse image of $\real$.
Choose $
\epsilon > 0 $ so that $ c_i + \epsilon < c_{i+1} - \epsilon $ for
$1 \leq i < r$.  Choose $A > 0$ so that $-A < c_1$ and $c_r < A$.  Since a level
set of $\tf$ corresponding to a regular value is a union of circles,
\begin{equation}
\label{eqn2}
\euler(M) = \euler( \{ \tf \leq -A \} \cup \{ \tf \geq A \} )
+ \sum_i \euler( \{ c_i - \epsilon \leq \tf \leq c_i + \epsilon \} )
\end{equation}
where $\euler$ denotes Euler characteristic.
The set $\{ \tf \leq -A \} \cup \{ \tf \geq A \}$ is homotopy
equivalent to the set $ \tf^{-1}(\infinity) $.  This is a connected
set, and is homotopy equivalent to a join of circles.  These circles
are the exceptional
sets where $\tf = \infinity$ together with the proper transform of
$\reallineinfinity$.  
Thus 
$$
\chi( \{ \tf \leq -A \} \cup \{ \tf \geq A \} ) = -\xi_{\rlinf,\infinity}(\tilde{f})
$$

Next, $M$ is a connected sum of copies of $\projective^2$, so 
$$\chi(M) = 1 - ( \xi_{\rlinf}(\tilde{f}) 
+ \xi_{\rlinf, \infinity}(\tf) + \xi_{\reallineinfinity,nc}(\tilde{f}) ) $$

At a critical value $c_i$, $\chi( \{ c_i - \epsilon \leq \tf \leq c_i
+ \epsilon \} )$ is the sum of the indices of the corresponding critical points, by Morse theory.
The sum of all these indices can be split into the parts
coming from critical points in the finite plane and the line at infinity.
Using this fact and the two equations above changes Equation (\ref{eqn2}) to
$$ 1 -  \xi_{\rlinf}(\tf) - \xi_{\rlinf, \infinity} (\tf) - \xi_{\reallineinfinity,nc}(\tilde{f}) 
= - \xi_{\rlinf,\infinity}(\tf) + i + i_{\reallineinfinity}(\tilde{f})$$
which proves (1).

\medskip

\noindent Proof of (2):
(See Figures \ref{std-crpt-l} and \ref{std-crpt-res-l}.)
Choose $\epsilon > 0$ so that $c$ is the only critical value in $(c-\epsilon,
c+\epsilon)$.
Let $C'$ be the (closed) exterior of the circle $C$ in the plane.
Let $N'$ be the connected component of  
$ \{(x,y) \in \real^2: c-\epsilon \leq f(x,y) \leq c+\epsilon \}
\cap C' $
containing $p$ in its closure.
Choose the circle $C$ large enough so that each boundary component of
$N'$ consists of
an arc of $f = c \pm \epsilon$ followed by an arc of $C$ followed by
an arc of $f = c \pm \epsilon$.

\begin{figure}
\postscript{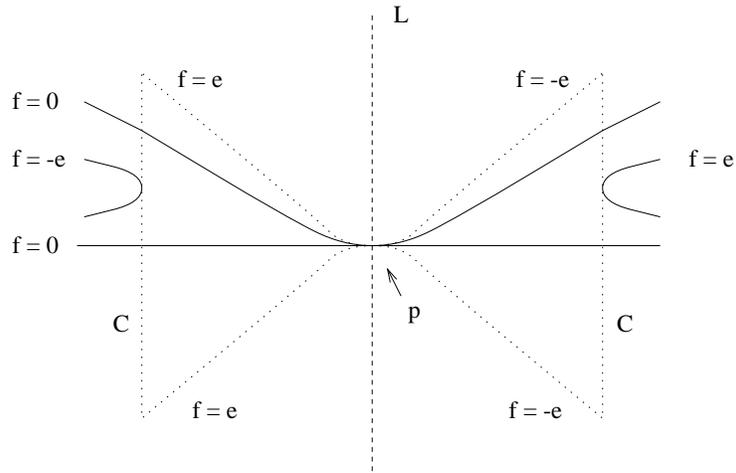}{0.7}
\caption{The region $N'$ 
(bounded by dotted lines) for the polynomial $y(xy-1)$ at $p = [1,0,0]$
and $c =0$.}
\label{std-crpt-l}
\end{figure}
\begin{figure}
\postscript{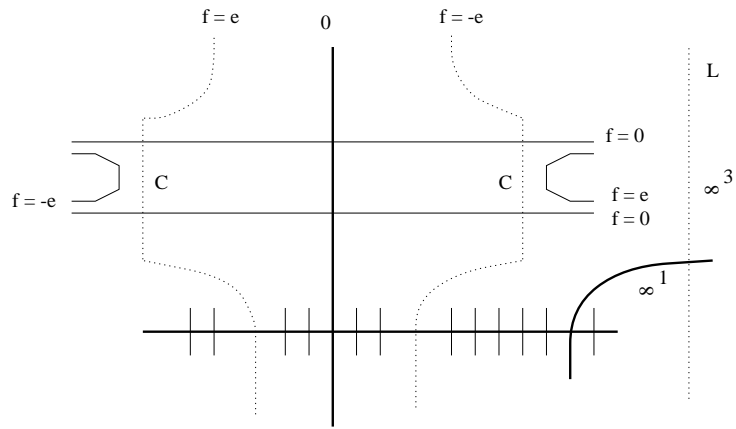}{0.7}
\caption{The region $N$ (bounded by dotted lines)}
\label{std-crpt-res-l}
\end{figure}

Let
$$N = \overline{\pi^{-1}(N')} \subset M$$
We assume that $N$ is connected; if it is disconnected the proof is similar.

We need a variant of the Poincar\'{e}-Hopf Theorem
for vector fields on a manifold, or more properly, a variant of Morse
theory on manifolds with boundary. 
(See, for instance, \cite{Milnor-TDV}, p. 35).
For an oriented manifold $X$ with boundary, the Euler characteristic
$\euler(X)$ is given by
$$\euler(X) = \sum \{\mbox{indices of internal critical points}\}+
\{\mbox{index on boundary}\} $$
where the index of the vector field on the boundary is measured with respect
to the outward pointing normal vector.
This result is true for a gradient vector field on a nonorientable
two-manifold $X$ without boundary, provided that the index is defined to be $+1$ at a local
extremum and $-1$ at a saddle, or more generally defined at an arbitrary
critical point using the result of Arnold \cite{Arnold-78} that the index
of a polynomial $f(x,y)$ at a point  $p$ is $1-r$, where $r$ is the number of real branches at $p$ of the curve $f(x,y) = f(p)$.
If $X$ has a boundary with an orientable collar neighborhood, then the
result is still true, provided that the index on the boundary is
measured according as Figure \ref{u-sign}.
Finally, the form we will use for $N$ is
$$\euler(N) = 1+ \sum \{\mbox{indices of internal critical points}\}+
\{\mbox{index on boundary}\} $$
The term $+1$ comes from the fact that $N$ has four corners (see
Figure \ref{std-crpt-res-l}).

Choose a Riemannian metric on $N$ so that it agrees, on the boundary components of $N$ consisting of arcs
of $C$, with the standard
metric on the plane.  We apply the above result to the gradient vector field of
$\tf$.
In the interior of $N$ there are the exceptional sets with $\tilde{f} = c$ 
and those critical points of $\tilde{f}$ which have critical value $c$.
Since $N$ retracts to the exceptional sets contained in it,
$$\euler(N) = 1 - \xi_{p,c}(\tilde{f})$$
The index of the internal critical points of $\tilde{f}$ is
$i_{p,c}(\tilde{f})$.
Finally, the index of the gradient vector field on the boundary
of $N$ is $i_{p,c}$. 
Combining these facts proves (2).

\medskip

\noindent Proof of (3):
(See Figure \ref{large-level-curves}.)
If $|f| = A$ is compact for $A \gg 0$, then Proposition
\ref{geom-prop-cpt} and Proposition \ref{index-prop} 
prove the result.
Hence we can suppose that $|f| = A$ is not compact.
If $\gamma$ is an end of the curve of tangencies and $c(\gamma) = 
\infinity$, then $\gamma$ intersects $|f| = A$ for
$A >> 0$.  

Let $I$ be a connected component of $|f| = A$ in $\real^2$.
Each $\gamma$ which meets $I$ has $c(\gamma) = \infinity$.
Since $I$ begins and ends outside $C$ by Proposition
\ref{geom-prop-noncpt} part (1), clearly
the sum of the $ i(\gamma)$ over all $\gamma$ meeting $I$ is $-1/2$.

By Proposition \ref{geom-prop-noncpt} part (2), 
there is a two-to-one correspondence between connected components
of $|f| = A$ and exceptional sets where $\tilde{f}$ is nonconstant.
This proves (3).
\end{xproof}

\begin{example}
Consider the resolution $\tf$ of the
polynomial $y(xy-1)$ shown in Figure \ref{std-crpt-res-g}.
Over $p = [1,0,0]$, $\tf$ has two
saddle points with critical value $0$, so $i_{[1,0,0],
0}(\tf) = -2$.
There is one exceptional set $E$ over $[1,0,0]$ with $\tf | E = 0$,
so $\xi_{[1,0,0], 0}(\tf) = 1$, and one exceptional set where $\tf$ is
nonconstant, so so $\xi_{[1,0,0], nc}(\tf) = 1$.
Over $p = [0,1,0]$, there is just one exceptional set where $\tf$ is
nonconstant, so $\xi_{[0,1,0], nc}(\tf) = 1$.
\end{example}

Recall that $f_d$ is the homogeneous term of highest degree of the
polynomial $f$, and that $\dr$ is the real degree of $f$ as defined in the
Introduction.

\begin{corollary}
\label{corollary-dr}
If $f_d$ has no repeated real linear factors, then $i = 1-\dr$.
Also, $i_{p,c}= 0$ for $p \in \rlinf$ and $c \in
\real$, and $i_{\rlinf, \infinity} = -\dr$.
\end{corollary}

\begin{figure}
\postscript{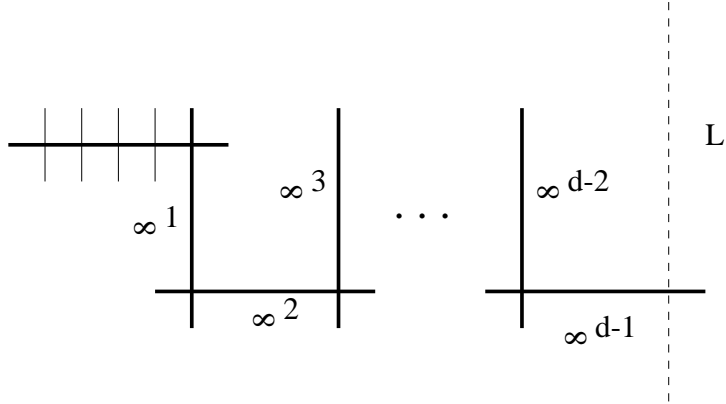}{0.7}
\caption{Resolution of a point where $d_p = 1$}
\label{linear-res-fig}
\end{figure}

\begin{xproof}
This ``geometrically obvious'' result follows from Proposition
\ref{res-index-formula}, since over each point where the level curves
of $f$ intersect $\rlinf$ the resolution is as in Figure
\ref{linear-res-fig}. 
\end{xproof}

Let $h(x,y)$ be a real polynomial whose homogeneous term of highest
degree is a product of irreducible real quadratic factors (ie, $\dr = 0$).
Instead of using a family of
concentric circles (the level curves of the polynomial $x^2+y^2$) to
define the $i_{p,c}$ , we could use the level curves of the polynomial
$h(x,y)$.  We let $i^h_{p,c}$ be the decomposition of the index
defined this way.

\begin{proposition}
\label{generic-h}
Let $f(x,y)$ be a real polynomial with isolated critical points.
If $h(x,y)$ is a real polynomial with $\dr = 0$, then $i^h_{p,c} = i_{p,c}$
for all $p \in \rlinf$ and $c \in \real$.
\end{proposition}

\begin{xproof}
Let $e$ be the degree of $h$.
We have that $i_{p,c} = i^{h'}_{p,c}$ for $h' = (x^2+y^2)^{e/2}$.

Choose a family $h^s$ of polynomials with $h^0 = h$ and $h^1 = h'$,
and so that $h^s$ is a polynomial in $x$, $y$ and $s$ with degree
$e$ in $x$ and $y$
and the homogeneous term of highest degree $(h^s)_e$ is a
product of irreducible real quadratic factors.
The curve of tangencies for $h^s$ is $f_xh^s_y
- f_yh^s_x=0$.
Choose a resolution $\tf : M \to \real$ of $f$, with $\pi: M \to
\projective^2$.
The curve of tangencies for $h^s$ lifts to $M$.

Let $q^s$ be 
an intersection point of the proper transform of the
lifted curve of tangencies with the exceptional set
$\pi^{-1}(\rlinf)$.
If
$\pi(q^s) = p$ and $\tf(q^s)=c$,
then $q^s$ contributes to
$i^{h^s}_{p,c}$.
The point $q^s$ varies continuously with $s$.
Clearly $\pi(q^s)$ is independent of $s$.
We will show that $\tf(q^s)$ is also independent of $s$.

If $q^s$ does not move as a function of $s$, then this is true.
Suppose $q^s$ moves.  Fix an $s$ and call it $s_0$.  Suppose that
$q^{s_0}$ is contained in an exceptional
set $E$. We will show that $\tf | E$ is constant. 
The function $\tf | E$ is rational.
By moving $s_0$ a little, we may assume that $q^{s_0}$
is not an intersection point of $\pi^{-1}(L)$ and that $q^{s_0}$ is not a critical
point of $\tf | E$.  For each $s$, the level sets of $\tilde{h^s}$
form a system of regular neighborhoods of $\pi^{-1}(L)$.  The level
curves of $\tilde{h^s}$ and $\tf$ are tangent along the curve of
tangencies for $h^s$.  Thus the level curves of $\tf$ become tangent
to $E$ in a neighborhood of $s_0$.  Hence $\tf$ is constant in a
neighborhood of $s_0$, and hence on $E$.
\end{xproof}

The proposition is obviously not true for $i_{p,\infinity}$ with 
$p \in \lineinf'$.
Lastly, is there a polynomial $f$ with a resolution $\tf$ and a point $q$ in
the exceptional set such that $\tf$ has a local extremum at $q$?


\section{A formula for $i$ in terms of a Morsification}


This section contains a formula (Proposition
\ref{morsification-formula}) for computing the index $i$ of a
polynomial and the
$i_{p,c}$ from Section 3 in terms of a Morsification.
The proof is straight-forward, and the results will not be used later.
Some examples and a conjecture are included.

\begin{definition}
A {\em deformation} of a real polynomial $f(x,y)$ of degree $d$ with
isolated critical points is a real
polynomial $h(x,y,s)$ of degree $d$ in $x$ and $y$ with $h(x,y,0) = f(x,y)$. 
We let $f^s(x,y) = h(x,y,s)$.
A deformation is a {\em
Morsification} if for small $s \neq 0$, $(f^s)_d$ 
(the homogeneous term of degree $d$ of $f^s$) has no repeated real
linear factors and the critical
points of $f^s$ in $\real^2$ are nondegenerate.
\end{definition}

It is easy to show that 
a polynomial of degree $d$ has a Morsification, and that the set of Morsifications is a dense open subset of the set of polynomials of degree $\leq d$.

If $f^s$ is a Morsification of $f$ and $p \in \rlinf$ is in the
closure of some level set of $f$ and $c \in \realinf$, we let
\begin{itemize}
\item $\tilde{d}_p(f^s)$ be the number of real linear factors in
$(f^s)_d$ which are deformations of the
factor corresponding to $p$ in $f_d$.  (The number $\tilde{d}_p(f^s)$ is also the number of points on
$\reallineinfinity$ through which the level sets of the
Morsification pass and which go to $p$ as $s \to 0$.)
\item $i^{\infinity}_{p,c}(f^s)$ be the
index of the critical points of $f^s$ which go to $p$ and whose
critical value goes to $c$ as $s \to 0$
\item $i^{\infinity}_{p}(f^s) 
= \sum_{c \in \realinf} i^{\infinity}_{p,c}(f^s)$
\item $i^{\infinity}(f^s) 
= \sum_{p \in \rlinf} i^{\infinity}_{p}(f^s)$
\end{itemize}

Note that these invariants depend not just on the expression for $f^s$
but also on the sign of $s$.

\begin{proposition}
\label{morsification-formula}
Let $f^s$ be a Morsification of $f$.
\begin{enumerate}
\item $i = 1- \dr(f^s) - i^{\infinity}(f^s)$
\item  $i_p = 1-\tilde{d}_p(f^s) -i^{\infinity}_p(f^s)$ for $p \in
\rlinf''$ (ie, for $p$ in the closure of some level set of $f$).
\end{enumerate}
\end{proposition}

\begin{xproof}
(1):  We have
$i(f^s) = i^{\infinity}(f^s)+i$ by summing the indices of the
critical points of $f^s$, and 
$i(f^s) = 1 - \dr(f^s)$ by Corollary \ref{corollary-dr}. 

\noindent (2):  
We may assume that $p = [1,0,0]$.
Choose a four-sided region on the $x > 0$ side
of the plane containing the critical points of $f^s$ on that side
which go to $p$ as $s \to 0$, and such that the left side of the region is
a segment of the circle $C$ containing the points of tangency $q$ which
approach $p$, the top and bottom are level sets of $f^s$,
and whose right side is a segment of a larger circle $C'$.
Orient the boundary of this region counter-clockwise.
Choose a similar region on the left side of the plane.
The index of $grad \, f^s$ about these two regions is clearly $1 - i_p
- \tilde{d}_p(f^s)$.
It is also the sum of the indices of the critical points in the
interior of the regions, which is $i^{\infinity}_p(f^s)$.
\end{xproof}


There is no obvious formula for $i_{p,c}$; see Conjecture
\ref{ipc-conjecture} at the end of this section.

\begin{example}
Let $f(x,y) = y^2-x$ and $p = [1,0,0]$.
Here $i_{p,c} =0$ for $c \in \realinf$.
The Morsification $f^s(x,y) = y^2 +sx^2 -x$ 
has a critical point at $(1/2s,0)$ with critical
value $-1/4s$.
If $s>0$ then $\tilde{d}_{p}(f^s) = 0$ and the critical point is
a mininum, so $i^{\infinity}_{p,0}(f^s) = 1$.
If $s<0$ then $\tilde{d}_{p}(f^s) = 2$ and the critical point is
a saddle, so $i^{\infinity}_{[p,0}(f^s) = -1$.
\end{example}

\begin{example}
Let $f(x,y) = y(xy-1)$ and $p = [1,0,0]$.
Here $i_{p,0} = 1$ and $i_{p,c} = 0$ for $c \neq 1$.
Define a Morsification by $f^s(x,y) = (y-sx)(xy-1)$.
(This deformation simply tilts the line in the zero locus of $f$.)
We have $\tilde{d}_{p}(f^s) = 2$.
For $s > 0$, $f^s$ has two real nondegenerate critical
points:

\medskip

\begin{tabular}{|l|l|l|}  \hline
critical point & type & critical value \\ \hline 
$(+ 1/ \sqrt{s}, +
\sqrt{s})$ & saddle & $0$ \\ \hline 
$(- 1/ \sqrt{s}, - \sqrt{s})$ & saddle & $0$ \\ \hline
\end{tabular}

\medskip  

\noindent Thus 
$i^{\infinity}_{p,0}(f^s) = -2$.
\end{example}


\begin{example}
Let $f(x,y) = x(y^2-1)=x(y+1)(y-1)$ and $p = [1,0,0]$.
Here $i_{p,\infinity} = -1$ and $i_{p,c} = 0$ for $c \in \real$.
We give two Morsifications.
The first is $f^s(x,y) =
x(y-sx+1))(y+sx-1)$.  
We have $\tilde{d}_p(f^s) = 2$.
The function $f^s$ has 
critical points:

\medskip
\begin{tabular}{|l|l|l|}  \hline
critical point & type & critical value \\ \hline 
$(0,\pm 1)$ & saddle & $0$ \\ \hline
$(1/s,0)$ & saddle & $0$ \\ \hline
$(1/(3s), 0)$ & minimum ($s>0$) & $-4/(27s)$ \\ 
              & maximum ($s<0$) &             \\ \hline
\end{tabular}
\medskip

\noindent Thus $i^{\infinity}_{p,0} = -1$ and 
$i^{\infinity}_{p,\infinity} = 1$.
The second Morsification is $f^s(x,y) = x(y^2 -1) + sx^3$.
For $s > 0$, $\tilde{d}_p(f^s) = 0$ and $f^s$ has critical points

\medskip
\begin{tabular}{|l|l|l|}  \hline
critical point & type & critical value \\ \hline 
$(0,\pm 1)$ & saddle & $0$ \\ \hline
$(1/\sqrt{3s},0)$ & mimimum & $-2/3\sqrt{3s}$ \\ \hline
$(-1/\sqrt{3s}, 0)$ & maximum  & $2/3\sqrt{3s}$ \\ \hline
\end{tabular}
\medskip

\noindent Thus $i^{\infinity}_{p,\infinity} = 2$.
(If $s < 0$, the only real critical points of $f^s$ are at $(0, \pm 1)$, and
$\tilde{d}_p(f^s) = 2$.) 
\end{example}



There is no obvious formula for $i_{p,c}$ in terms of a deformation,
as can be seen in the case of $x(y^2-1)$ above.
However, the following conjecture seem reasonable:

\begin{conjecture}
\label{ipc-conjecture}
If $f^s$ is a deformation of $f$ and $p \in \rlinf''$, then
$i_{p,c} \leq -i^{\infinity}_{p,c}(f^s)$ for $c \in \real$, and
$i_{p,\infinity} \leq 1 - i^{\infinity}_{p,\infinity}(f^s)$.
\end{conjecture}

For a deformation $f^s$ of $f$ it is easy to find bounds on the number
of local maxima, minima and saddles near a point $p \in \rlinf$.  It
would be interesting to see what possible combinations of these can
occur, similar to the investigation in \cite{REU} or \cite{Shustin-96}.


\section{Bounds on $i$}

This section contains the main results of this paper, the bounds on
the index $i$ of the gradient vector field of a real polynomial.
The main tool is a bound on $i_p$ (Lemma \ref{i-d-local-estimate}).
This, together with the interpretation of $i_{\rlinf, \infinity}$ in
terms of a resolution (Proposition \ref{res-index-formula}) 
and some lemmas using techniques from Section 4 give
the first main result (Theorem \ref{dr-theorem}).
We next give a refinement (Lemma \ref{lemma-2}) of Lemma \ref{i-d-local-estimate}. 
This and a number of technical details gives the second 
main result (Theorem \ref{max-theorem}).
As remarked in the Introduction, these is still a large gap between
the upper bounds and known examples.

Let $f(x,y)$ be a real polynomial of degree $d$ with isolated critical points.
For $p \in \rlinf$ and $c \in \realinf$, recall that 
$$i_{p,c} = \sum i(\gamma)$$ 
where the sum is over ends $\gamma$ of the curve of tangencies with $p(\gamma) = p$ and
$c(\gamma) = c$.
We let
$$i_{p,c}^{abs} = \sum |i(\gamma)| $$
and
$$i^{abs}_p = \sum_{c \in \realinf} i^{abs}_{p,c}$$
These invariants can be computed from a resolution of $f$, and in
particular are integers, although
this is not evident from the proof of \ref{res-index-formula}.  
(The invariants $i_{p,c}^{N}$ and so forth later in this section can
also be computed from a resolution.)

Recall that $\lineinf'$ and $\lineinf''$ were introduced in Section 4,
that $\lineinf'' \subset \lineinf'$ 
and that $l''$ is the number of elements in $\lineinf''$.
If $p = [a,b,0] \in \rlinf$, we let $d_p$ be the
multiplicity of the factor $(bx-ay)$ in $f_d$.
The next result is the local analogue of the estimate $|i| \leq d-1$ from the
Introduction.  
This estimate, like the global one, is proved by relating the index to
an algebraic intersection number.

\begin{lemma}
\label{i-d-local-estimate}
If $p \in \lineinf'$ (ie, if $d_p > 0$), then 
$$i_p^{abs} \leq  d_p - 1 $$
\end{lemma}

This follows from the next two results.
We let $\Gamma$ be the projective completion of the curve of
tangencies, and $\Gamma_\complex$ be its complexification.
We use $(A,B)_p$ to denote the intersection number of the curves $A$
and $B$ at $p$.

\begin{lemma}
\label{i-d-local-estimate-1}
$i^{abs}_p \leq (\Gamma_\complex \intersect \clineinf)_p$
\end{lemma}

\begin{xproof}
The number $i^{abs}_p$ is at most one half the number of ends
$\gamma$ of the curve of tangencies with $p(\gamma) = p$.
This number is the number of real branches of the completion of the
curve of tangencies at $p$,
which is at most the number of branches of
$\Gamma_\complex$ at $p$.
This number is at most $(\Gamma_\complex \intersect \clineinf)_p$, 
since no component of the curve of
tangencies is contained in $\clineinf$.
\end{xproof}

\begin{lemma}
\label{i-d-local-estimate-2}
$(\Gamma_\complex \intersect \clineinf)_p = d_p -1$
\end{lemma}

\begin{xproof}
Without loss of generality we may assume that $p = [1,0,0]$.
We have that $f = y^{d_p}h(x,y)
+{\mbox { \{terms of lower order\} } }$ where $d_p \geq 1$, $h(x,y)$ is homogeneous of
degree $d-d_p$, and
$y$ does not divide
$h(x,y)$.
Changing coordinates to $p$ and computing as
in \cite[III.3]{Fulton} shows that 
$(\Gamma_\complex \intersect
\clineinf)_p = d_p - 1$
\end{xproof}

Lemma \ref{i-d-local-estimate} is sharp:
The polynomial of Example \ref{many-max-min-ex} at $p = [1,0,0]$ has $d_p = k$ and 
$i^{abs}_{p} = i_{p,0} = k-1$.
 Another example is provided by the
polynomial 
$x(y+1)(y+2) \dots (y+k)$
at $p = [1,0,0]$, which has
$d_p = k$ and $i^{abs}_{p} = -i_{p,0} = k-1$.

Recall that the real degree of the polynomial is 
$$\dr = \sum_{p \in \lineinf} d_p $$ 
We let
$$\drt = \sum_{p \in \lineinf''} d_p $$ 
This is the sum of the $d_p$'s over those $p$ in the line at infinity
which are in the closure of some level set of the polynomial $f$.
Note that
\begin{equation}
\label{string}
l'' \leq \drt \leq \dr \leq d
\end{equation}

\begin{proposition}
\label{i-upper-bound}
If $f(x,y)$ is a real polynomial with isolated critical points, then
$$i \leq 1 + \drt - 2l''$$
\end{proposition}

\begin{xproof}
Let $\tf$ be a resolution of $f$.
We have 
$$i = 1 + \sum_{\scriptstyle p \in \rlinf \atop \scriptstyle 
c \in \real} i_{p,c} + i_{\reallineinfinity, \infinity}$$
$$\leq  1 + \sum_{\scriptstyle p \in \rlinf'' \atop \scriptstyle 
c \in \real} i_{p,c}^{abs} - \xi_{\reallineinfinity, nc}(\tf)$$
$$\leq 1 +\sum_{p \in \rlinf''} (d_p -1 ) - \xi_{\reallineinfinity, nc}(\tf)$$
$$ \leq  1 + \drt - l'' - \xi_{\reallineinfinity, nc}(\tf)$$
The first line follows from Proposition 
\ref{geometric-index-formula}.
The second follows since 
$i_{\rlinf, \infinity} = - \xi_{\reallineinfinity, nc}(\tf)$ 
by Part (3) of Proposition \ref{res-index-formula}, 
$i_{p,c} \leq i_{p,c}^{abs}$ by definition, and
$i_{p,c}^{abs} = 0$ for $p \in \rlinf - \rlinf''$ by Lemma \ref{lemma-i-p-c}.
The third line follows from  
Lemma \ref{i-d-local-estimate}. 
The result follows from Corollary \ref{l-xi}.
\end{xproof}


To get a lower bound on the index, we need to compactify the plane
$\real^2$ by the circle 
$$ \circleinf = \{ (a,b,0) \in (\real^3 - 0) / \real^+ \} $$ 
The projection map 
$$\circleinf \to \lineinf$$
which takes $(a,b,0)$ to $[a,b,0]$ will be denoted by $q \mapsto |q|$.
If $\gamma$ is a end of the curve of tangencies, we let $q(\gamma) \in
\circleinf$ be the endpoint of the closure of $\gamma$ in
$\circleinf$.  

For example, for $y-(xy-1)^2$ there are ends $\gamma$ and $\gamma'$
with $q(\gamma) = (1,0,0)$, $q(\gamma') = (-1,0,0)$, $c(\gamma) =
c(\gamma') = 0$ and $i(\gamma) = + 1/2$, $i(\gamma') = - 1/2$.
The two cancel out to give $i_{[1,0,0],0} = 0$.

Let
$$\circleinf'' = \{ q \in \circleinf : \mbox{ \ There is a $t \in \real$ such that
$q$ is in the closure of $f(x,y) = t$} \} $$
This is a finite set of points; we let $s''$ denote the number of
points in this set.
Since the fibers of the projection map
$$\circleinf'' \to \lineinf''$$
consist of one or two points, 
$$l'' \leq s'' \leq 2l''$$
Thus the string of inequalities (\ref{string}) becomes
\begin{equation}
\label{string-2}
0 \leq \frac{1}{2} l'' \leq \frac{1}{2} s'' \leq l'' \leq \drt \leq \dr \leq d
\end{equation}

The next two lemmas are preparation for proving
Proposition \ref{i-lower-bound}.

\begin{lemma}
$$\biggl | \sum_{q(\gamma) \in \circleinf''} i(\gamma)\biggr | \leq 
\drt - l''$$
\end{lemma}

\begin{xproof}
$$\biggl | \sum_{q(\gamma) \in \circleinf''} i(\gamma)
\biggr | \leq 
\sum_{q(\gamma) \in \circleinf''} |i(\gamma)| $$
$$\leq \sum_{p \in \lineinf''} \sum_{| q(\gamma)| = p} |i(\gamma)| $$
$$\leq \sum_{p \in \lineinf''} (d_p - 1)$$
$$= \drt - l''$$
\end{xproof}

\begin{lemma}
$$ \biggl| \sum_{q(\gamma) \in \circleinf - \circleinf''} i(\gamma)
\biggr| \leq \frac{1}{2} s''$$
\end{lemma}

\begin{xproof}
If $\circleinf''$ is empty ($s'' = 0$), then all the level sets of $f$
are compact and the left-hand sum is zero by Proposition
\ref{index-prop}.

Now suppose that $\circleinf''$ is not empty.
Fix a connected component $V$ of $\circleinf - \circleinf''$.
Let $A \gg 0$ and
let $I(V)$ be the connected component of $|f| = A$ which goes to $V$ as
$A$ goes to infinity.
By Propositions \ref{geom-prop-cpt} and \ref{geom-prop-noncpt}, $I(V)$ is
not compact and hence has its endpoints on $\circleinf$.
We have (see the proof of Part (3) of Proposition
\ref{res-index-formula}.)
that
$$\sum i(\gamma) = - \frac{1}{2}$$
where the sum is over all ends $\gamma$ of the curve of tangencies
which intersect $I(V)$.
If an end $\gamma$ has $q(\gamma) \in V$, then $\gamma$ intersects $I(V)$, and
no other connected component of $|f| = A$.
(Since the closure of no level curves of $f$ pass through the
endpoint of $\gamma$, the function $f$ approaches infinity
monotonely on
$\gamma$. (Lemma \ref{lemma-cgamma}))
However, some of ends $\gamma$ of the curve of tangencies may
intersect $I(V)$ but have their endpoints on the endpoints of $V$.
(I know of no examples of polynomials with this property, though.)
The sum of the $i(\gamma)$ which intersect $I(V)$ and have endpoint a chosen
endpoint of $V$ is 0 or $+ 1/2$.
The sum over both endpoints of $V$ is thus 0, $+1/2$ or 1.
Hence
$$\sum_{q(\gamma) \in V} i(\gamma) = - \frac{1}{2},\ 0 \mbox{ or } +\frac{1}{2} $$

Since $\circleinf - \circleinf''$ has $s''$ connected components, this
proves the lemma.
\end{xproof}


\begin{proposition}
\label{i-lower-bound}
If $f(x,y)$ is a real polynomial with isolated critical points, then
$$i \geq 1 - \dr$$
\end{proposition}

\begin{xproof}
We have 
$$i = 1 + \sum_{\gamma} i(\gamma)$$
$$= 1 + \sum_{q(\gamma) \in \circleinf''} i(\gamma)
+ \sum_{q(\gamma) \in \circleinf - \circleinf''} i(\gamma)
$$
$$\geq 1 + l'' - \drt -\frac{1}{2}s''$$
$$\geq 1 - \drt$$
$$\geq 1 - \dr$$
The third line follows from the two lemmas above, and the fourth and
fifth by the string of inequalities (\ref{string-2}).
\end{xproof}

The following is our first main result.

\begin{theorem}
\label{dr-theorem}
Let $f(x,y)$ be a real polynomial of real degree $\dr$ with isolated
critical points, and let $i$ be the index of $grad \, f$ around a large
circle containing the critical points.
If all the level sets of $f$ are compact, then $i = 1$.  Otherwise
$$|i| \leq \dr - 1$$
\end{theorem}

\begin{xproof}
If the level sets are compact then the result is obvious (see
Proposition \ref{index-prop}.)
If some level sets are not compact, then $l'' >0$.
The upper bound follows from Proposition \ref{i-upper-bound}, and the
lower bound from Proposition \ref{i-lower-bound}.
\end{xproof}

As remarked in the Introduction and Corollary
\ref{corollary-dr}, it is easy to find ``generic'' polynomials
which realize the lower bound of this theorem.
The upper bound appears too high; the estimate of Proposition
\ref{i-upper-bound} is somewhat better.
Finally, the result seems somewhat obvious and one could hope for a
better proof.

We now further decompose $i_{p,c}$ and its refinements defined above.
For $p \in \rlinf$ and $c \in \realinf$, 
recall that $i_{p,c} = \sum i(\gamma)$, summed over all
ends $\gamma$ of the curve of tangencies with $p(\gamma) = p$ and
$c(\gamma) = c$.
We let $i^T_{p,c}$
(respectively, $i^N_{p,c}$) be the sum of the $i(\gamma)$'s 
such that the corresponding curve $\gamma$ is tangent
(respectively, not tangent) to $\rlinf$ at $p$.
Thus
$$i_{p,c} = i^N_{p,c} + i^T_{p,c} $$
We similarly decompose $i_{p,c}^{abs}$.
As before, these numbers are all integers.
For example, the polynomial $y(xy-1)$ of Example \ref{std-crpt-ex} has $i_{[1,0,0],0} =
i^{N}_{[1,0,0],0} = 1$, and
the polynomial in Example \ref{two-parabola-ex} has $i_{[1,0,0],0} =
i^{T}_{[1,0,0],0} = 1$.
We also let
$$i^{N, abs}_p = \sum_{c \in \realinf} i^{N, abs}_{p,c}$$
and define $i^{T, abs}_p$ similarly.
The following lemma is a refinement of Lemma \ref{i-d-local-estimate}.

\begin{lemma}
\label{lemma-2}
If $p \in \rlinf'$, then
$$i^{N,abs}_{p} 
+ 2 i^{T,abs}_{p}   \leq d_p-1$$
\end{lemma}

\begin{xproof}
We let $\Gamma^T$ (respectively, $\Gamma^N$) be the product of the branches of
the curve of tangencies $\Gamma$ tangent (respectively, not tangent) to $\rlinf$ at $p$, so that
$\Gamma = \Gamma^T \Gamma^N$ near $p$.
As in Lemma \ref{i-d-local-estimate-1}, we have that
\begin{equation}
\label{eqn-n}
i^{N,abs}_p \leq (\Gamma^N_\complex \intersect \clineinf )_p
\end{equation} 
Similarly
\begin{equation}
\label{eqn-t}
i^{T,abs}_p \leq \frac{1}{2}(\Gamma^T_\complex \intersect \clineinf
)_p
\end{equation}
since these branches are
tangent to $\lineinf$ at $p$.
Thus 
$$i^{N,abs}_p + 2i^{T,abs}_p \leq (\Gamma^N_\complex \intersect
 \clineinf )_p + (\Gamma^T_\complex \intersect \clineinf )_p = (\Gamma_\complex \intersect
 \clineinf )_p = d_p -1$$
as before.
\end{xproof}

Some stronger local estimates are probably true.
In fact, let $f = f^Nf^T$ at $p$, where $f^T$ (respectively, $f^N$) are
the branches of $f = t$ tangent (respectively, not tangent) to $\rlinf$ at
$p$ (which is independent of $t \in \real$), and let
$d_p^T$ (respectively, $d_p^N$) be the intersection number of 
$f^T = t$ (respectively, $f^N = t$) with $\rlinf$ at $p$.
Thus
$$d_p = d_p^N + d_p^T$$
It seems reasonable to expect that 
$i^{N,abs}_{p} \leq d^N_p - 1$
and
$i^{T,abs}_{p} \leq (1/2) d^T_p - 1$
for $p \in \rlinf$, and that these estimates are sharp.

We need one more technical lemma:

\begin{lemma}
\label{lemma-n-t}
Fix $p \in \rlinf$ and let $\tf$ be a resolution of $f$.
If there are $c, c' \in \real$ such that $i^{N,abs}_{p,c} > 0$ and
$i^{T,abs}_{p,c'}>0$,
then $\xi_{p,nc}(\tf) \geq 2$
(ie, there are at least two exceptional sets over $p$ on which $\tf$
is not constant).
\end{lemma}

\begin{xproof}
There are ends $\gamma$ and $\gamma'$ of the curve of
tangencies with $p(\gamma) = p(\gamma') = p$,
$c(\gamma) = c$ and $c(\gamma') = c'$, and 
with $\gamma$ (respectively $\gamma'$) tangent (respectively, not tangent) to $\rlinf$ at $p$.
By the proof of Lemma \ref{lemma-i-p-c}, the closure of the level
curve $f = c$ (respectively, $f = c'$) intersects an exceptional set
$E$ (respectively, $E'$) over $p$.
By the proof of Lemma \ref{geom-lemma-1}, there is at least one
exceptional set over $p$ where $\tf$ is not constant.
In fact, there are at least two such exceptional sets: 
Since 
$\gamma$ and $\gamma'$ 
have distinct tangents at $p$, 
the limit of $f$ must be infinite on all but a finite
number of tangent directions between these by \cite[Proposition
1.3]{REU}, so $E$ and $E'$ are distinct and the chain of exceptional sets connecting $E$ and $E'$
must have at least one member $E_0$ with $f | E_0 = \infinity$.
Thus there must be an exceptional set in the chain connecting $E$ and
$E_0$ where $\tf$ is nonconstant, and similarly between $E'$ and $E_0$.
\end{xproof}

\begin{theorem}
\label{max-theorem}
If $f(x,y)$ is a real polynomial of degree $d$ with isolated
critical points, and $i$ is the index of $grad \, f$ around a large
circle containing the critical points, then
$$ i \leq max \{ 1,d-3 \}$$
\end{theorem}

The difficult part of this proof is the case when $l''=0$, i.e. when the
closures of the level curves of the polynomial intersect the line at
infinity at just one point.

\begin{xproof}
If $l'' = 0$ then $i = 1$ by Lemma \ref{index-prop}.
If $l'' \geq 2$ then $i \leq \dr -3$ by Proposition \ref{i-upper-bound}.
Thus we must treat the case $l'' = 1$.
We may suppose without loss of generality that $p = [1,0,0]$ is the
only point where the real level curves of $f$ intersect $\rlinf$.
The point $p$ will remain fixed for the rest of the proof.

Let $\tf$ be a resolution of $f$.
From Proposition \ref{geometric-index-formula}, Lemma \ref{lemma-i-p-c} and part (3) of Proposition
\ref{res-index-formula} we have that
\begin{equation}
\label{equation-l-1}
i = 1 + \sum_{c \in \real} i_{p,c} - \xi_{p,nc}(\tf)
\end{equation}
Since $\xi_{\rlinf, nc}(\tf) \geq l'' = 1$ by Corollary \ref{l-xi},
we also have the weaker form of this equation:
$$i \leq \sum_{c \in \real} i_{p,c}$$

Suppose $d_p < d$.  We have that $d_p \leq d-2$ since the roots of $f_d$ other than
$p$ are
complex and hence occur in conjugate pairs.  Thus 
$$i \leq \sum_{c \in \real} i_{p,c} \leq 
\sum_{c \in \real} i^{abs}_{p,c} \leq i^{abs}_p
\leq d_p -1 \leq d-3$$
where the fourth inequality follows from Lemma
\ref{i-d-local-estimate}.

Thus we may assume that $d_p = d$,
so that
$$f(x,y) = \pm y^d + h(x,y)$$ 
where $h$ has degree $e < d$.
If $h$ is a function of $x$ alone, then 
from far away $f(x,y)$ looks like $\pm y^d \pm x^e$, which has $i = 0$
or $\pm 1$.
Thus we may assume that $h$ is a nonconstant function of both $x$ and
$y$.

The rest of the proof is divided into three cases:

\noindent Case 1:  Suppose $\sum_{c \in \real} i^N_{p,c} = 0$. 
Then  
$$i \leq \sum_{c \in \real} i_{p,c} = \sum_{c \in \real} i^T_{p,c}
\leq \sum_{c \in \real} i^{T,abs}_{p,c}
\leq i^{T,abs}_p$$
$$\leq \frac{1}{2}(d_p-1)
\leq max \, \{1, d-3 \}
$$
where the fourth inequality follows from Lemma \ref{lemma-2}.

\noindent Case 2:  Suppose that
$\sum_{c \in \real} i_{p,c}^T = 0$. 
We have that
$$i \leq \sum_{c \in \real} i_{p,c}
= \sum_{c \in \real} i_{p,c}^N 
\leq \sum_{c \in \real} i_{p,c}^{N,abs}
\leq i^{N,abs}_p
\leq ( \Gamma^N_\complex  \intersect \clineinf ) _p$$ 
where the last inequality is Equation (\ref{eqn-n}).

Since 
$f(x,y) = \pm y^d + h(x,y)$ 
where $h$ is a nonconstant function of both $x$ and $y$ of degree less than $d$,
a computation shows that $z$ divides the term of lowest degree
in the curve of tangencies localized at $p$.
Hence $\Gamma^T$ is nonempty, so 
$( \Gamma^T_\complex \intersect \clineinf ) _p \geq 2$.

Thus
$$( \Gamma^N_\complex  \intersect \clineinf ) _p =
( \Gamma_\complex  \intersect \clineinf  )_p -
( \Gamma^T_\complex  \intersect \clineinf  )_p$$
$$ \leq (d_p - 1) - 2 $$
$$ = d_p - 3$$
$$ = d-3$$

\noindent Case 3:
Suppose that $i^{N, abs}_{p,c} > 0$ and 
$i^{T, abs}_{p,c'}>0$ for some $c, c' \in \real$.
We have by Equation (\ref{equation-l-1}) and Lemma \ref{lemma-n-t}
that
$$i  \leq \sum_{c \in \real} i_{p,c} - 1$$
$$ \leq \sum_{c \in \real} i_{p,c}^{abs} - 1$$
$$ = \sum_{c \in \real} i_{p,c}^{N,abs} + 
2 \biggl( \sum_{c \in \real} i_{p,c}^{T,abs} \biggr) 
- \sum_{c \in \real} i_{p,c}^{T,abs} - 1$$
$$ \leq (d-1) -1-1 = d-3 $$
where the last inequality follows from Lemma \ref{lemma-2}.
\end{xproof}




\section{Vanishing cycles}

Suppose $p \in \rlinf$ 
and $c \in \real \cup \{ \infinity \}$.
In this section we relate the term $i_{p,c}$ in the decomposition of
the index $i$ of a real polynomial $f(x,y)$ to the number of vanishing
cycles $\nu_{p,c}$
of the corresponding complex polynomial at $(p,c)$.
This number
is defined to be the ``jump'' in the Milnor number of the family
$f(x,y)=t$ of complex polynomials at $p$ when $t=c$.
(For a detailed discussion of this notion, in particular the case $c =
\infinity$, the reader is referred to \cite{Durfee-P96}.)
Recall that $i^{abs}_{p,c}$ was defined in the last
section, and that $i_{p,c} \leq i^{abs}_{p,c}$.

\begin{proposition}
Let $f(x,y)$ be a real polynomial with isolated critical points.
If $p \in \rlinf'$ (ie, $p$ is a zero of $f_d$)  and $c
\in \realinf$, then
$$i_{p,c}^{abs} \leq \nu_{p,c} $$
\end{proposition}

\begin{xproof}
Suppose $c \in \real$.
(The proof for $c = \infinity$ is similar.)
As in the proof of Lemma \ref{i-d-local-estimate-1}, 
The number
$i^{abs}_{p,c}$ is at most one half the number of ends
$\gamma$ of the curve of tangencies with $p(\gamma) = p$ and
$c(\gamma)= c$.
Since $f$ is either strictly increasing or
decreasing on each end $\gamma$ (Lemma
\ref{lemma-cgamma}),
we may assume without loss of generality (replace $f$ by $-f$), that
the number of ends $\gamma$ with $p(\gamma) = p$,
$c(\gamma)= c$ and $f|\gamma < c$ is at most 
the number of ends $\gamma$ with $p(\gamma) = p$,
$c(\gamma)= c$ and $f|\gamma > c$.
Hence $i^{abs}_{p,c}$ is at most  
the number $v$ of ends $\gamma$ with $p(\gamma) = p$,
$c(\gamma)= c$ and $f|\gamma > c$.

Since $f$ is strictly decreasing to $c$ along $\gamma$, $v$ is the
number of intersection
points in $\real^2$ of the 
curves $xf_y - yf_x = 0$ (the curve of tangencies) and $f = c +
\epsilon$ which approach $p$ as $\epsilon \downarrow 0$.
If we assume without loss of generality that $p = [1,0,0]$,
we may replace the curve $xf_y - yf_x = 0$ by the curve $f_y = 0$.   
Thus $v$ is at most the number of intersection points of the
complex curves $f_y = 0$ and $f = c + \epsilon$ which approach $p$ as
$\epsilon \downarrow 0$. 
This number
is well-known to be $\nu_{p,c}$ 
(see for example \cite[2.13]{Durfee-P96}).
\end{xproof}

The inequality of the proposition can be strict; for example the
polynomial $y(x^ay-1)$ at $p = [1,0,0]$ and $c=0$ has $i_{p,c} = 1$
and $\nu_{p,c} = a+1$.  

\bibliographystyle{alpha}

\newcommand{\etalchar}[1]{$^{#1}$}

\bigskip

\noindent Department of Mathematics

\noindent Mount Holyoke College

\noindent South Hadley, MA 01075

\medskip

\noindent email:  adurfee@mtholyoke.edu

\end{document}